\documentclass[10pt]{amsart}
\input{diagrams}

\newtheorem{proposition}{Proposition}[section]
\newtheorem{lemma}[proposition]{Lemma}
\newtheorem{corollary}[proposition]{Corollary}
\newtheorem{theorem}[proposition]{Theorem}

\theoremstyle{definition}
\newtheorem{definition}[proposition]{Definition}
\newtheorem{example}[proposition]{Example}

\theoremstyle{remark}

\def \rhd{\hspace*{-1mm}\rightharpoondown\hspace*{-1mm}}

\newcommand{\thlabel}[1]{\label{th:#1}}
\newcommand{\thref}[1]{Theorem~\ref{th:#1}}
\newcommand{\selabel}[1]{\label{se:#1}}
\newcommand{\seref}[1]{Section~\ref{se:#1}}
\newcommand{\lelabel}[1]{\label{le:#1}}
\newcommand{\leref}[1]{Lemma~\ref{le:#1}}
\newcommand{\prlabel}[1]{\label{pr:#1}}
\newcommand{\prref}[1]{Proposition~\ref{pr:#1}}
\newcommand{\colabel}[1]{\label{co:#1}}

\newcommand{\exlabel}[1]{\label{ex:#1}}
\newcommand{\exref}[1]{Example~\ref{ex:#1}}
\newcommand{\delabel}[1]{\label{de:#1}}
\newcommand{\deref}[1]{Definition~\ref{de:#1}}
\newcommand{\eqlabel}[1]{\label{eq:#1}}

\newcommand{\Hom}{{\rm Hom}}
\newcommand{\End}{{\rm End}}

\newcommand{\Ker}{{\rm Ker}\,}

\def\ot{\otimes}

\def\Br{{\rm Br}}
\def\BM{{\rm BM}}
\def\BT{{\rm Br}'}
\def\BTM{{\rm BM}'}

\def\Mm{{\mathcal M}}

\def\End{{\rm End}}
\def\Gal{{\rm Gal}}
\def\Br{{\rm Br}}
\def\BM{{\rm BM}}

\def\Hom{{\rm Hom}}
\def\E{{\rm E}}
\def\Pic{{\rm Pic}}

\def\text#1{{\rm {\rm #1}}}

\def\ol{\overline}
\def\ul{\underline}

\def\*C{{{}^*\mathcal C}}

\begin{document}
\title[The equivariant Brauer group of a group]{The equivariant Brauer group of a group}
\author{S. Caenepeel}
\address{Faculty of Applied Sciences,
Vrije Universiteit Brussel, VUB, B-1050 Brussels, Belgium}
\email{scaenepe@vub.ac.be}
\urladdr{http://homepages.vub.ac.be/\textrm{\~}scaenepe/}
\author{F. Van Oystaeyen}
\address{Department of Mathematics and Computer Science\\ University of Antwerp.
 Middelheimlaan 1, B-2020 Antwerp, Belgium}
\email{fred.vanoystaeyen@ua.ac.be}
\author{Y.H. Zhang}
\address{School of Mathematics and Computing Science,
Victoria  University of Wellington,
Wellington, New Zealand}
\email{yinhuo.zhang@vuw.ac.nz}

\begin{abstract}
We consider the Brauer group ${\BTM}(k,G)$ of a
group $G$ (finite or infinite) over a commutative ring $k$ with
identity. A split exact sequence
$$1\longrightarrow  \BT(k)\longrightarrow \BTM(k,G)\longrightarrow\Gal(k,G)
\longrightarrow 1$$ is obtained. This generalizes the Fr\"ohlich-Wall
exact  sequence (\cite{F,FW1})from the case of a field to the case of a commutative ring, and generalizes the Picco-Platzeck exact sequence (\cite{PP}) from the finite case of $G$ to the infinite case of $G$. Here $\BT(k)$ is the Brauer-Taylor 
group of  Azumaya algebras (not necessarily with unit) cf. \cite{Ta}. The method developed in this paper might provide a key to computing the equivariant Brauer group of an infinite quantum group introduced in \cite{CVZ}.
\end{abstract}
\maketitle

\section*{Introduction}
In \cite{FW1, FW2}, Fr\"ohlich and Wall introduced the {\it equivariant
Brauer group} $\BM(k,G)$ of a finite group $G$ and a commutative ring $k$,
consisting of equivalence classes of $k$-Azumaya algebras on which $G$ acts as
a group of automorphisms. In \cite{F} Fr\"ohlich studied further the equivariant 
Brauer group of an arbitrary group over a field.  
$\BM(k,G)$ can be completely described by the
Brauer group $\Br(k)$ and the second cohomology group $H^2(G,U(k))$ via the following split exact sequence when $\Pic(k)=1$:
$$1\to \Br(k)\to \BM(k,G)\to H^2(G,U(k))\to 1.$$
If $\Pic(k)\neq 1$, then we replace $H^2(G,U(k))$ by the group
$\Gal(k,G)$ of $kG$-Galois objects (these are $G$-strongly graded rings having $k$
as part of degree $e$). This sequence is originally due to Picco and Platzeck \cite{PP},
and several variations and generalizations of it have appeared in the literature.
Beattie \cite{Be} generalized the sequence to the Brauer group of algebras with an
action by a finitely generated projective cocommutative Hopf algebra, and
Childs \cite{Ch} gave an exact sequence involving algebras
with a finite abelian group grading inducing a coaction, using a bilinear map.
All the generalizations (see \cite{Ca} for a detailed discussion) have in common that
the (co)acting group or Hopf algebra is finite.\\
The aim of this note is to extend the Picco-Platzeck exact sequence to infinite groups.
We can still introduce the invariants $\Gal(k,G)$, $H^2(G,U(k))$ and $\BM(k,G)$, and define a map
$\pi:\ \BM(k,G)\to \Gal(k,G)$, but the problem is that $\pi$ is in general not surjective if $G$ is infinite. The idea is then to replace the Brauer group by the Brauer-Taylor group. Taylor \cite{Ta} defined the {\it Brauer-Taylor group} $\BT(k)$, 
consisting of equivalence classes of Azumaya
algebras without a unit. If $k$ is noetherian, then there is an isomorphism between
$\BT(k)$ and the second {\'e}tale cohomology group. It is known that the
second {\'e}tale cohomology group, and, a fortiori,
Taylor's Brauer group is not always torsion, and it follows from 
Gabber's Theorem \cite{Ga,KO2} that the classical Brauer group $\Br(k)$ is exactly the torsion part of $\BT(k)$.\\
Let $G$ be a group, and consider Taylor-Azumaya algebras on which $G$ acts as a group of automorphisms. Such a
Taylor-Azumaya algebra is called a $G$-module Taylor-Azumaya algebra.
The set of isomorphism classes of $G$-module Taylor-Azumaya algebras modulo Morita
equivalence is a group, called the equivariant Brauer-Taylor group of $G$ and $k$.

In Sections \ref{se:1.1} and \ref{se:1.2}, we recall some facts about the Brauer-Taylor group and about the group of Galois objects. Our methods will be based on the theory of
multiplier algebras and multiplier
Hopf algebras, and these are briefly discussed in  Sections \ref{se:1.3} and \ref{se:1.4}.
The equivariant Brauer-Taylor group is defined in \seref{2}.
In \seref{3}, we define the group homomorphism $\pi$
from the Taylor-Brauer group to the group of $kG$-Galois objects and show
that $\pi$ is surjective and split. The main result is the exact
sequence in \thref{3.9}, describing the equivariant Brauer-Taylor group as
the direct sum of the Brauer-Taylor group and the group of $kG$-Galois objects.

We will use the following notation: if $C$ is an object of a category, then
$C$ will also be a notation for the identity morphism on $C$.

\section{Preliminaries}\selabel{1}
\subsection{Taylor-Azumaya algebras}\selabel{1.1}
In \cite{Ta}, Taylor introduced a Brauer group consisting of equivalence
classes of Azumaya algebras without a unit. This new invariant was 
investigated further in \cite{SG,RT}. Our main reference will be \cite{Ca}.\\
Let $k$ be a commutative ring, and $A$ a $k$-algebra, not necessarily with
unit. We call $A$ a {\sl unital} $k$-algebra if the canonical morphism
$A\ot_A A\to A$ is an isomorphism. In a similar way, a left $A$-module $M$
is called unital if the canonical map $A\ot_A M\to M$ is an isomorphism.
An algebra with unit is unital, and then any module is unital. If $A$
is a unital algebra, then ${}_A\Mm^u$ will be the category of unital left
$A$-modules. Unital right modules and bimodules are defined in a similar way.\\
If $A$ is a unital $k$-algebra, then the enveloping algebra $A^{\rm e}=
A\ot A^{\rm op}$ is also unital. The center of $A$ is defined as
$Z(A)=\End_{A^{\rm e}}(A)$.
If $A$ has a unit, then we recover the classical definition of the center.
From the fact that $A^2=A$, it follows that $Z(A)$ is a commutative
$k$-algebra with unit, and $A$ is a $Z(A)$-algebra. $A_l$ (resp. $A_r$)
will be our notation for $A$ when viewed as a left
(resp. a right) $A^{{\rm e}}$-module.
We can define (strict) Morita contexts connecting unital algebras in the
obvious way: we require that the connecting bimodules are unital bimodules
(see \cite[Sec. 1.1]{Ca}).

\begin{definition}\delabel{1.1}
Let $A$ be a unital, faithful $k$-algebra. $A$ is called a Taylor-Azumaya
algebra if it satisfies the following equivalent conditions.
\begin{enumerate}
\item There exists an invertible $k$-module $I$ such that the functors
$$F:\ {}_k\Mm\to {}_{A^{{\rm e}}}\Mm^u,~~F(N)=A_l \ot N$$ 
and
$$H:\ {}_{A^{{\rm e}}}\Mm^u\to{}_k\Mm,~~H(M)=(A_r\ot I)\ot{}_{A^{{\rm e}}} M$$
form a pair of inverse equivalences;
\item the pair of adjoint functors $(F,H_l)$, with
$$H_l:\ {}_{A^{{\rm e}}}\Mm^u\to{}_k\Mm,~~H_l(M)=\Hom_{A^{\rm e}}(A_l,M)$$
is a pair of inverse equivalences.
\end{enumerate}
\end{definition}

For the proof of the equivalence of the conditions, we refer to
\cite[Prop. 2.2.5 and 2.2.8]{Ca}. From the second condition, it follows immediately that a Taylor-Azumaya algebra is central. A Taylor-Azumaya algebra with a unit is
a usual Azumaya algebra. More conditions characterizing Taylor-Azumaya algebras
may be found in \cite[Sec. 2.2, 4.2 and 4.3]{Ca}. If $A$ is a Taylor-Azumaya
algebra, then the functors $H$ and $H_l$ from the definition are isomorphic.\\
Consider two $k$-modules $M$ and $M'$ and a surjective $k$-linear map
$\mu:\ M'\ot M\to k$. We call $\ul{M}=(M,M',\mu)$ a dual pair of
$k$-modules, and define the associated {\it elementary algebra}
$\E(\ul{M})=M\ot M'$, with multiplication
\begin{equation}\eqlabel{1.1.1}
(m\ot m')(n\ot n')=\mu(m'\ot n)m\ot n' ,
\end{equation}
for $m,n\in M$ and $m',n'\in M'$.
Then $\E(\ul{M})$ is a Taylor-Azumaya algebra. $M$ is a unital left $\E(\ul{M})$-module,
and $M'$ is a unital right $\E(\ul{M})$-module; the actions are given by the
formulas
\begin{equation}\eqlabel{1.1.2}
(m\ot m')n=m\mu(m' \ot n)~~{\rm and}~~n'(m\ot m')=\mu(n'\ot m)m'.
\end{equation}
If $\E(\ul{M})$ has a unit, then
$M$ and $M'$ are faithfully projective as $k$-modules, $M'\cong M^*$, 
and $\mu:\ M'\ot M\cong M^*\ot M\to k$ is the evaluation map. In this case,
$\E(\ul{M})\cong \End_k(M)$.

\begin{proposition}\prlabel{1.2} \cite[Prop. 3.1.1]{Ca}
Let $A$ and $B$ be Taylor-Azumaya algebras. The following assertions are
equivalent.
\begin{enumerate}
\item $A$ and $B$ are Morita equivalent;
\item $A\ot B^{\rm op}$ is an elementary algebra;
\item there exist dual pairs $\ul{M}$ and $\ul N$ such that
$A\ot \E(\ul{M})\cong B\ot \E(\ul{N})$ as $k$-algebras.
\end{enumerate}
\end{proposition}

Now consider the set $\BT(k)$ of Morita equivalence classes 
represented by a {\sl finitely generated}
Taylor-Azumaya algebra. $\BT(k)$ is a group under the operation induced by the
tensor product (see \cite[Sec. 3.2]{Ca}). We recall from \cite{SG} that every
{\sl flat} Taylor-Azumaya algebra is Morita equivalent to a finitely generated
Taylor-Azumaya algebra. There exist examples of Taylor-Azumaya algebras that
are not flat (see \cite{SG} or \cite[Sec. 4.3]{Ca}). If $k$ is a field, then
$\BT(k)=\Br(k)$.

\subsection{Strongly graded rings}\selabel{1.2}
Let $k$ be a commutative ring, and $G$ an arbitrary group. A $G$-graded $k$-algebra
$S$ is called a $kG$-{\sl Galois object} if the canonical map
$$\gamma:\ S\ot S\to S\ot kG,~~\gamma(s\ot t)=\sum_{g\in G} st_g\ot g$$
is an isomorphism. This is equivalent to (see \cite[Prop. 8.2.1]{Ca})
\begin{itemize}
\item $S$ is strongly graded, that is, $S_gS_h= S_{gh}$, for all $g,h\in G$;
\item $S_e=k$, where $e$ is the unit element of $G$.
\end{itemize}
The set of equivalence classes $\Gal(k,G)$
of $kG$-Galois objects forms a group under the
operation induced by the cotensor product:
$$S\square T=\bigoplus_{g\in G} S_g\ot T_g.$$

\subsection{Multipliers}\selabel{1.3}
From \cite{VD}, we recall the notion of multiplier algebra of an algebra $A$.
A multiplier of $A$ is a pair
$(\rho_1, \rho_2)$ consisting of a right $A$-linear map $\rho_1$ and a left
$A$-linear map $\rho_2$ such that $a \rho_1(b) =
\rho_2 (a)b$ for all $a,b \in A$. If
 $x = (\rho_1, \rho_2)$ is a multiplier, then we write $xb = \rho_1(b)$ and $bx = \rho_2
(b)$. $M(A)$, the set of all  multipliers on $A$, is a $k$-algebra with
multiplication $(\rho_1,\rho_2)(\rho'_1,\rho'_2)=(\rho_1\circ \rho'_1,
\rho'_2\circ \rho_2)$ and unit $(A,A)$. We have an algebra map
$A\to M(A)$; the multiplier corresponding to $a\in A$ is given by left and right
multiplication. 
$M(A)$ is an algebra with unit $(A,A)$. If $A$ has identity, then $A=M(A)$.
Indeed, if $x=(\rho_1,\rho_2)\in M(A)$, then $\rho_1(1)=1\rho_1(1)=\rho_2(1)1=
\rho_2(1)$, and $\rho_1(b)=1\rho_1(b)=\rho_2(1)b$, $\rho_2(c)=c\rho_1(1)$.
Let us compute the multiplier algebra of an elementary algebra.

\begin{proposition}\prlabel{1.3}
Let $\ul{M}=(M,M',\mu)$ be a dual pair, and $E=E(\ul{M})$ the corresponding
elementary algebra. Then the multiplier algebra $M(E)$ is isomorphic to
$$\ol{E}=\{(f,f')\in E_1\times E_2~\mid~ \mu(m'\ot f(m))=\mu(f'(m')\ot m),
~{\rm for~all~} m\in M, m'\in M'\},$$
where $E_1=\End(M)$ and $E_2=\End(M')^{op}$.
\end{proposition}

\begin{proof}
We have a map
$$\alpha:\ \ol{E}\to M(E),~~\alpha(f,f')=x=(\rho_1=f\ot M',\rho_2=M\ot f').$$
It is straightforward to see that $x$ is a multiplier.\\
Next we define
$$\beta:\ M(E)\to \ol{E},~~\beta(\rho_1,\rho_2)=(f,f'),$$
with $f:\ M\to M$ and $f':\ M'\to M'$ are defined as follows:
$$f(m)=\sum_i \rho_1(m\ot m'_i)m_i~~{\rm and}~~
f'(m')=\sum_i m'_i\rho_2(m_i\ot m'),$$
where $\mu(\sum_i m'_i\ot m_i)=1$. We compute that
$$\rho_1(m\ot m')=\sum_i \rho_1((m\ot m'_i)(m_i\ot m'))=
\sum_i \rho_1(m\ot m'_i)m_i\ot m',$$
$$\rho_2(m\ot m')=\sum_i \rho_2((m\ot m'_i)(m_i\ot m'))=
\sum_i m\ot m'_i\rho_1(m_i\ot m').$$
Since $(\rho_1,\rho_2)$ is a multiplier, we have, for all $m,n\in M$
and $m',n'\in M'$, that
\begin{eqnarray*}
&&\hspace*{-2cm}
(m\ot m')\rho_1(n\ot n')=\mu(m'\ot f(n))m\ot n'\\
&=& \rho_2(m\ot m')n\ot n'=\mu(f'(m)\ot n)m\ot n',
\end{eqnarray*}
and it follows from the fact that $E$ is faithful that $(f,f')\in \ol{E}$.
finally, it is easy to see that $\alpha$ is the inverse of $\beta$.
\end{proof}

\subsection{Multiplier Hopf algebras}\selabel{1.4}
From \cite{VD}, we recall the definition of a multiplier Hopf algebra. If $A$ is a unital algebra, then $A\ot A$ is again a unital algebra, with non-degenerate product. The tensor product $M(A)\ot M(A)$ is in a natural way a subalgebra of $M(A\ot A)$, and the embedding $A\to M(A\ot A)$ factors through this map. A {\sl multiplier Hopf algebra} is a unital algebra $A$ together with an algebra map $\Delta:~A\ot A\to M(A\ot A)$ such that $\Delta(a) (1 \otimes b)$ and $(a \otimes 1)\Delta(b)$ belong to $A\otimes A$ for all $a,b \in A$ and  $\Delta$ is coassociative in the sense that
$$(a \otimes 1 \otimes 1)(\Delta \otimes A)(\Delta(b)(1 \otimes c)) =(A \otimes \Delta)((a \otimes 1) \Delta (b))(1 \otimes 1 \otimes c)$$
for all $a,b,c\in A$. In this paper, we will need the following example of multiplier Hopf algebra.
\begin{example}\exlabel{1.4}
Let $G$ be a group. Let $(kG)^*$ be the dual of the group algebra $kG$.
Let $p_g:\ kG\to k$ be the map defined by $p_g(h)=\delta_{g,h}$,
and $k(G)=\oplus_{g\in G} kp_g$. Then $M(k(G))=(kG)^*$.
In a similar way, $M(k(G\times G))=(kG\times kG)^*$. There is a natural identification $k(G)\ot k(G)\cong k(G\times G)$, and $M(k(G)\ot k(G))\cong M(k(G\times G))$. 
The multiplication on $G$ induces an algebra map $\Delta~: k(G)\to M(k(G) \otimes k(G))$ given by $(\Delta(f))(p,q) = f(pq)$. One may use (informal) Sweedler sigma notation for the comultiplication:  
$$\Delta(p_g)=\sum {p_g}_{(1)}\ot {p_g}_{(2)}=\sum_{h\in G}p_h\ot p_{h^{-1}g}$$ 
for $p_g\in k(G)$.  Moreover, we have that 
$\Delta (f_1)(1 \otimes f_2) \in k(G) \otimes k(G)$ and $(f_1 \otimes 1)\Delta(f_2) \in k(G) \otimes k(G)$, for all $f_1,f_2 \in k(G)$. 
From the associativity of the multiplication on $G$, it easily follows that $\Delta$ is coassociative. Thus $(k(G),\Delta)$ is a multiplier Hopf algebra.
\end{example}

\subsection{The equivariant Brauer group}\selabel{1.5}

Let $A$ be a $k$-algebra and $G$ be a group. If $G$ acts as a group of $k$-algebra automorphisms on $A$ with trivial $G$-action on $k$, then
we say that $A$ is a $G$-module algebra. A Taylor-Azumaya algebra $A$ that is
also a $G$-module algebra is called a $G$-module Taylor-Azumaya algebra. If
$\ul{M}=(M,M',\mu)$ is a dual pair, such that $M$ and $M'$ are $G$-modules, and
$$\mu(g\cdot m'\ot g\cdot m)=\mu(m'\ot m),$$
for all $g\in G$, $m\in M$ and $m'\in M'$, then $\E(\ul{M})$ is a 
$G$-module Taylor-Azumaya algebra, with diagonal action $g\cdot (m\ot m')=
g\cdot m\ot g\cdot m'$. We say that $\E(\ul{M})$ is an elementary
$G$-module algebra.\\
Let $A$ and $B$ be unital $G$-module algebras. We can define (strict)
Morita contexts connecting $A$ and $B$, requiring that the bimodules in
the Morita contexts are $G$-modules that are unital on both sides, and
that the connecting bimodule maps preserve the $G$-action. We have the
following generalization of \prref{1.2}. The proof is an easy adaption of
the proof of \cite[Prop. 3.1.1]{Ca}.

\begin{proposition}\prlabel{2.1}
Let $A$ and $B$ be $G$-module Taylor-Azumaya algebras. The following assertions are
equivalent.
\begin{enumerate}
\item $A$ and $B$ are Morita equivalent as $G$-module algebras;
\item $A\ot B^{\rm op}$ is an elementary $G$-module algebra;
\item there exist dual pairs of $G$-modules $\ul{M}$ and $\ul M$ such that
$A\ot \E(\ul{M})\cong B\ot \E(\ul{N})$ as $G$-module algebras.
\end{enumerate}
\end{proposition}

Now we introduce the equivariant Brauer group $\BTM(k,G)$ as the set of
Morita equivalence classes of $G$-module Taylor-Azumaya algebras, represented
by a finitely generated $G$-module Taylor-Azumaya algebra. It is a group
under the operation induced by the tensor product. If $k$ is a field, then
$\BTM(k,G)\cong \BM(k,G)$, the Brauer group of $G$-module Azumaya algebras.

\section{The split exact sequence}\selabel{3}
Let $k$ be a commutative ring,
$A$ a unital $k$-algebra, and $G$ a group. ${}_k\Mm^G$ will be the category of
graded $k$-modules, and ${}_A\Mm^{uG}$ the category of $G$-graded unital
left $A$-modules. Take an invertible $k$-module $I$. The functors $F$ and $H$
from \deref{1.1} send graded modules to graded modules. If $A$ is a Taylor-Azumaya
algebra, then $H\cong H_l$, hence $H_l$ also sends graded modules to graded
modules, and we have a pair of inverse equivalences $(F,H_l)$, with
$$F:\ {}_k\Mm^G\to {}_{A^{{\rm e}}}\Mm^{uG},~~F(N)=A_l \ot N$$ 
and
$$H_l:\ {}_{A^{{\rm e}}}\Mm^{uG}\to{}_k\Mm^G,~~H_l(M)=\Hom_{A^{\rm e}}(A,M).$$
In particular, if $M$ is a $G$-graded unital left $A^{\rm e}$-module, then
$\Hom_{A^{\rm e}}(A,M)$ is also $G$-graded. This means that every left
$A^{\rm e}$-linear map $f:\ A\to M$ can be decomposed into a finite sum
$f=\sum_{g\in G} f_g$, with $f_g(A)\subset M_g$, for all $g$.\\

Let $A$ be a $G$-module algebra, and consider the smash product
$S_A=A\# kG$, with multiplication
$$(a\# g)(b\# h)=a(g\cdot b)\# gh,$$for $a, b\in A$ and $g, h\in G$. We have an algebra map $\eta:\ A\to S_A$,
$\eta(a)=a\# 1_G$, making $A$ into a subalgebra of $S_A$. If $A$ is a unital algebra, then $S_A$ is a $G$-graded
unital left $A^{\rm e}$-module, with left action $(a\ot a')\cdot (b\# g)=
ab(g\cdot a')\# g$, and with $G$-grading given by ${\rm deg}(b\# g)=g$.

\begin{lemma}\lelabel{3.1}
Let $A$ be a $G$-module Taylor-Azumaya algebra. Then $\Hom_{A^{\rm e}}(A,S_A)=
H_l(S_A)$ is a $kG$-Galois object.
\end{lemma}

\begin{proof}
We have already seen at the beginning of this section that $H_l(M)$ is a
$G$-graded module, for every $G$-graded unital left $A^{\rm e}$-module $M$.
Take $\varphi,\psi\in \Hom_{A^{\rm e}}(A,S_A)$. The product of $f$ and $g$ is defined
as follows: for all $a\in A$, take $a_i,a'_i\in A$ such that $a=\sum_i a_ia'_i$,
and put
$$(\varphi\psi)(a)= \sum_i \varphi(a_i)\psi(a'_i).$$
The product is well-defined because $A\cong A\ot_A A$, and both $\phi$ and $\psi$ are $A$-bimodule maps. It is easy to see that $\eta$ is a unit for this multiplication, 
and that
$\Hom_{A^{\rm e}}(A,S_A)$ is a $G$-graded algebra.
It remains to be shown that the canonical map
$$\gamma:\ \Hom_{A^{\rm e}}(A,S_A)\ot \Hom_{A^{\rm e}}(A,S_A)\to
\Hom_{A^{\rm e}}(A,S_A)\ot kG,~~\gamma(\varphi\ot \psi)=\varphi\psi\ot h,$$
if $\psi$ is homogeneous of degree $h$. Since $(F,H_l)$ is an equivalence
of categories, it suffices to show that $F(\gamma)=A\ot \gamma$ is an
isomorphism. Since $(F,H_l)$ is an equivalence of categories, the evaluation
map
$${\rm ev}:\ A\ot  \Hom_{A^{\rm e}}(A,S_A)\to S_A,~~{\rm ev}(a\ot\varphi)= \varphi(a)$$
is an isomorphism. Hence we have isomorphisms
\begin{eqnarray*}
&&\hspace*{-2cm}
A\ot  \Hom_{A^{\rm e}}(A,S_A)\ot  \Hom_{A^{\rm e}}(A,S_A)\cong
S_A\ot  \Hom_{A^{\rm e}}(A,S_A)\\
&\cong&
S_A\ot_A A\ot  \Hom_{A^{\rm e}}(A,S_A)\cong S_A\ot_A S_A.
\end{eqnarray*}
The composition of these isomorphisms is called $\alpha$, and is given by
$$\alpha(ab\ot \varphi\ot\psi)=\varphi(a)\ot_A\psi(b).$$
Also 
$${\rm ev}\ot kG:\ A\ot  \Hom_{A^{\rm e}}(A,S_A)\ot kG\to S_A\ot kG$$
is an isomorphism, and we define $\beta$ by the commutativity of the
following diagram:
$$\begin{diagram}
A\ot\Hom_{A^{\rm e}}(A,S_A)\ot \Hom_{A^{\rm e}}(A,S_A)&\rTo^{A\ot \gamma}&
A\ot\Hom_{A^{\rm e}}(A,S_A)\ot kG\\
\dTo^{\alpha}&&\dTo_{{\rm ev}\ot kG}\\
S_A\ot_A S_A&\rTo^{\beta}& S_A\ot kG
\end{diagram}$$
It now suffices to show that $\beta$ is an isomorphism. Take $a\# g,
b\# h\in S_A$, and write ${\rm ev}^{-1}(a\# g)=\sum_i a_i\ot \varphi_i$,
${\rm ev}^{-1}(b\# h)=\sum_j b_j\ot \psi_j$, with $\varphi_i$ and 
$\psi_j$ homogeneous of respectively degree $g$ and $h$. Then
$$\alpha(\sum_{i,j} a_ib_j\ot \varphi_i\ot \psi_j)=
\sum_{i,j} \varphi_i(a_i)\ot \psi_j(b_j)=(a\# g)\ot_A(b\#h),$$
hence
\begin{eqnarray*}
&&\hspace*{-2cm}
\beta((a\# g)\ot_A(b\#h))=
({\rm ev}\ot kG)(A\ot \gamma)(\sum_{i,j} a_ib_j\ot \varphi_i\ot \psi_j)\\
&=& ({\rm ev}\ot kG)(\sum_{i,j} a_ib_j\ot \varphi_i \psi_j\ot h)
= \sum_{i,j} (\varphi_i \psi_j)(a_ib_j)\ot  h\\
&=& (\sum_i\varphi_i(a_i))(\sum_j\varphi_j(b_j))
= (a\# g)(b\#h) =a (g\cdot b)\# gh.
\end{eqnarray*}
A straightforward computation shows that the inverse of $\beta$ is given by
the formula
$$\beta^{-1}((ab\# g)\ot h)=(a\# gh^{-1})\ot_A (hg^{-1})\cdot b\# h.$$
\end{proof}

We will now present some alternative descriptions of
$\Hom_{A^{\rm e}}(A,A\# kG)$. The first one is based on  the fact
that $\Hom_{A^{\rm e}}(A,A\# kG)$ is graded, and the second one 
is in terms of multipliers. Then we will show that there is an anti-algebra
homomorphism to the multipliers of $A$.

\begin{lemma}\lelabel{3.1a}
Let $A$ be a $G$-module Taylor-Azumaya algebra.
The part of degree $g$ of $\Hom_{A^{\rm e}}(A,A\# kG)$ is isomorphic to
$$\End_A^g(A)=\{f_g\in \End_{A-}(A)~|~f_g(ab)=f_g(a)(g\cdot b),~{\rm for~all}~a,b\in A\}.$$
Consequently, every $f\in \Hom_{A^{\rm e}}(A,A\# kG)$ can be written uniquely as
$$f(a)=\sum_{g\in G} f_g(a)\#g,$$
with $f_g\in \End_A^g(A)$, and only finitely many of the $f_g$ different from zero.
\end{lemma}

\begin{proof}
Assume that $f\in \Hom_{A^{\rm e}}(A,A\# kG)$ is homogeneous of degree $g$.
Then $f(a)\in A\# kg$, for every $a\in A$, so we can write $f(a)= f_g(a)\# g$.
From the fact that $f$ is left $A$-linear, it follows that $f_g$ is left $A$-linear.
From the fact that $f$ is right $A$-linear, we obtain
$$f_g(a)(g\cdot b)\# g= f(a)b=f(ab)=f_g(ab)\#g,$$
and the result follows.
\end{proof}

\begin{lemma}\lelabel{3.2}
Let $A$ be a $G$-module Taylor-Azumaya algebra. Then we have an
isomorphism of algebras
$$M(A\# kG)^A\cong \Hom_{A^{\rm e}}(A, A\# kG).$$
\end{lemma}

\begin{proof}
Take a multiplier $x=(\rho_1,\rho_2)\in M(A\# kG)$. Observe that
$x\in M(A\# kG)^A$ if and only if $\rho_1$ is left $A$-linear and
$\rho_2$ is right $A$-linear. We now define a map
$$\alpha:\ M(A\# kG)^A\to \Hom_{A^{\rm e}}(A, A\# kG),~~
\alpha(\rho_1,\rho_2)=f,$$
with $f(a)=\rho_2(a\#1)$. $f$ is left $A$-linear since $(\rho_1,\rho_2)$ 
is a multiplier,
and from the fact that $\rho_2$ is also right $A$-linear, it follows that $f$
is right $A$-linear. Thus $f$ is an $A$-bimodule map and hence a left $A^e$-map. 
Next we define
$$\beta:\ \Hom_{A^{\rm e}}(A, A\# kG)\to M(A\# kG)^A,~~
\beta(f)=(\rho_1,\rho_2),$$
with
$$\rho_1(a\# g)=f(a)(1\# g)~~{\rm and}~~\rho_2(a\# g)=(1\# g)f(g^{-1}\cdot a).$$
We remark that $1\# g\in M(A\# kG)$. $\rho_1$ is right $A\# kG$-linear since
\begin{eqnarray*}
&&\hspace*{-2cm}
\rho_1((a\# g)(b\# h))= \rho_1(a(g\cdot b)\# gh)= f(a(g\cdot b))(1\# gh)\\
&=& f(a)(g\cdot b\# gh)= \rho_1(a\# g)(b\# h).
\end{eqnarray*}
In a similar fashion, we can show that $\rho_2$ is left $A$-linear.
We next compute that
\begin{eqnarray*}
&&\hspace*{-15mm} (a\# g)\rho_1(b\# h)= (a\# g)f(b)(1\# h)=
(1\# g)(g^{-1}\cdot a\# 1)f(b) (1\# h)\\
&=& (1\# g)f((g^{-1}\cdot a)b) (1\# h)= (1\# g)f(g^{-1}\cdot a) (b\# h)=
\rho_2(a\# g)(b\# h),
\end{eqnarray*}
and it follows that $(\rho_1,\rho_2)$ is a multiplier. From the fact that
$f$ is left $A$-linear, it follows immediately that $\rho_1$ is left $A$-linear.
$\rho_2$ is right $A$-linear since
\begin{eqnarray*}
&&\hspace*{-2cm} \rho_2((a\# g)b)=\rho_2(a(g\cdot b)\# g)=
(1\# g) f((g^{-1}\cdot a)b)\\
&=& (1\# g) f(g^{-1}\cdot a)b=\rho_2(a\# g)b,
\end{eqnarray*}
where we now used the fact that $f$ is right $A$-linear. This shows that
$\beta(f)=(\rho_1,\rho_2)\in M(A\# kG)^A$. It is clear that
$\alpha$ is the inverse of $\beta$. Finally, $\alpha$ is an algebra homomorphism.
Indeed, if $\alpha(\rho_1,\rho_2)=f$ and $\alpha(\rho'_1,\rho'_2)=f'$,
then
\begin{eqnarray*}
&&\hspace*{-15mm}
\alpha(\rho_1\circ \rho'_1,\rho'_2\circ\rho_2)(ab)=
\rho'_2(\rho_2(ab\# 1))= \rho'_2(\rho_2(a\# 1)b)\\
&=&\rho'_2(\rho_2(a\# 1)(b\# 1))
= \rho_2(a\# 1)\rho'_2(b\# 1)=f(a)f'(b)=(ff')(ab).
\end{eqnarray*}
where we used the fact that $\rho_2$ is right $A$-linear and $\rho'_2$
is left $A\#kG$-linear.
\end{proof}

\begin{lemma}\lelabel{3.2a}
We have an algebra anti-homomorphism
$$p:\Hom_{A^{\rm e}}(A, A\# kG)\to M(A).$$
\end{lemma}

\begin{proof}
Take $f\in \Hom_{A^{\rm e}}(A, A\# kG)$, and define $p(f)=(\sigma_1,\sigma_2)$
by
$$\sigma_1(a)=\sum_{h\in G}  f_h(h^{-1}\cdot a)~~{\rm and}~~
\sigma_2(a)=\sum_{h\in G}  f_h(a),$$
where we use the notation introduced in \leref{3.1a}. It is straightforward to
show that $(\sigma_1,\sigma_2)$ is a multiplier. Let us show that
$p$ is an anti-homomorphism. Take $f,f'\in \Hom_{A^{\rm e}}(A, A\# kG)$,
and let $p(f)=(\sigma_1,\sigma_2)$, $p(f')=(\sigma'_1,\sigma'_2)$,
$p(ff')=(\tau_1,\tau_2)$. We have to show that
$\tau_2=\sigma_2\circ \sigma'_2~~{\rm and}~~\tau_1=\sigma'_1\circ \sigma_1$.
It suffices to prove this for $f$ and $f'$ homogeneous, respectively of
degree $g$ and $h$. We then have that
$f(a)=f_g(a)\# g~~{\rm and}~~f'(a)=f'_h(a)\# h$,
and we compute
\begin{eqnarray*}
&&\hspace*{-2cm}
(ff')(ab)=(f_g(a)\# g)(f'_h(b)\#h)=f_g(a)(g\cdot f'_h(b))\# gh\\
&=& f_g(af'_h(b))\# gh= f_g(f'_h(ab))\# gh.
\end{eqnarray*}
It follows immediately that
$$\tau_2(ab)=f_g(f'_h(ab))= (\sigma_2\circ \sigma'_2)(ab).$$
We finally compute
\begin{eqnarray*}
&&\hspace*{-15mm}
(\sigma'_1\circ \sigma_1)(ab)= f'_h(h^{-1}\cdot f_g(g^{-1}\cdot(ab)))=
f'_h\Bigl(h^{-1}\cdot \bigl((g^{-1}\cdot a)f_g(g^{-1}\cdot b)\bigr)\Bigr)\\
&=& f'_h\Bigl(\bigl((h^{-1}g^{-1})\cdot a\bigr)
\bigl(h^{-1}\cdot f_g(g^{-1}\cdot b)\bigr)\Bigr)=
 f'_h\bigl((gh)^{-1}\cdot a\bigr) f_g(g^{-1}\cdot b)\\
&=& f_g\Bigl(f'_h\bigl((gh)^{-1}\cdot a\bigr)(g^{-1}\cdot b)\Bigr)=
f_g\bigl(f'_h((gh)^{-1}\cdot (ab))\bigr)=\tau_1(ab).
\end{eqnarray*}
\end{proof}

\begin{example}\exlabel{3.3}
Let $A$ be a Taylor-Azumaya algebra, with trivial $G$-action: $g\cdot a=a$,
for all $g\in G$ and $a\in A$. Then $A$ is a $G$-module Taylor-Azumaya algebra,
and $A\# kG= A\ot kG$ as $k$-algebras. Furthermore
$$\Hom_{A^{\rm e}}(A,A\ot kG)\cong \Hom_{A^{\rm e}\ot k}(A\ot k,A\ot kG)
\cong \Hom_{A^{\rm e}}(A,A)\ot kG\cong kG.$$
\end{example}

Recall from \cite{BCM} that an action by a group $G$ on an algebra $A$ with identity is called {\it strongly inner} if there is a group morphism $f:\ G\to U(A)$, the group of units of $A$, such that $g\cdot a=f(g)af(g^{-1})$ for all $g\in G$ and$a\in A$. Obviously, this definition makes no sense if $A$ has no unit. We then have the following generalization.

\begin{definition}\delabel{3.4}
Let $G$ be a group and $A$ a unital algebra. An action of $G$ on $A$ is called strongly inner if there exists a group morphism $f:\ G\to U(M(A))$, the group of units of $M(A)$, such that $g\cdot a=f(g)af(g^{-1})$, for all $g\in G$ and $a\in A$.
\end{definition}

It is easy to see that a strongly inner action of $G$ on an algebra $A$ can be extended to a strongly inner action on the multiplier algebra $M(A)$ of $A$.

\begin{lemma}\lelabel{3.5}
Let $\ul{M}=(M,M',\mu)$ be a dual pair of $G$-modules. Then the induced action of $G$ on $E(\ul{M})$ is strongly inner. Conversely, 
a strongly inner $G$-action on an elementary algebra $E(\ul{M})$  is induced by a $G$-action on $\ul{M}$.
\end{lemma}

\begin{proof}
Let $\psi:\ G\to \End(M)$ and $\psi':\ G \to \End(M')$ be the representation maps of
$G$ on $M$ and $M'$. Then $(\psi(g),\psi'(g^{-1}))\in \ol{E}$ (see \prref{1.3})
since
$$\mu(m\ot g^{-1}\cdot m')=\mu(g\cdot m\ot (gg^{-1})\cdot m')=\mu(g\cdot m\ot m').$$
Let $f(g)$ be the  multiplier corresponding to $(\psi(g),\psi'(g^{-1}))$, as in
\prref{1.3}. Then we compute
$$g\cdot (m\ot m')=\psi(g)(m)\ot \psi'(g)(m')=f(g)(m\ot m')f(g^{-1}).$$
Conversely, assume that $G$ acts strongly innerly on an elementary algebra $E$.
Then there exists a group morphism $\lambda:\ G\to U(M(E))\cong U(\ol{E})$, 
$\lambda(g)=(f(g),f'(g))\in \ol{E}$ such that
$$g\cdot (m\ot m')=\lambda(g)(m\ot m')\lambda(g^{-1})=f_g(m)\ot f'_{g^{-1}}(m').$$
Now define actions of $G$ on $M$ and $M'$ by $g\cdot m=f_g(m)$ and
$g\cdot m'=f'_g(m')$. These actions are associative because $F$ is a group
homomorphism. Using the fact that $(f(g),f'(g))\in \ol{E}$, we compute
$$\mu(g\cdot m'\ot g\cdot m)=\mu(f'_{g^{-1}}(m')\ot f_g(m))=
\mu(m'\ot f'_{g^{-1}}(f_g(m)))=\mu(m'\ot m),$$
so $\ul{M}$ is a dual pair of $G$-modules.
\end{proof}

\begin{lemma}\lelabel{3.6}
Let $A$ be a $G$-module Taylor-Azumaya algebra. The $G$-action on $A$ is
strongly inner if and only if $\Hom_{A^{\rm e}}(A,A\# kG)\cong kG$
as $G$-graded algebras.
\end{lemma}

\begin{proof}
First suppose that the $G$-action is strongly inner. There exists
$f:\ G\to M(A)$ such that $g\cdot a=f(g)af(g^{-1})$, for all
$g\in G$ and $a\in A$. The map
$$\alpha:\ A\ot kG\to A\# kG,~~\alpha(a\ot g)=af(g^{-1})\# g$$
is an isomorphism of $G$-graded unital $A^e$-modules. Consequently
$H_l(A\ot kG)\cong H_l(A\# kG)$ as $G$-graded modules, and it is easy to
see that this is a morphism of $G$-graded algebras. Now we have seen in
\exref{3.3} that $H_l(A\ot kG)=kG$, and one implication follows.\\
Conversely, assume that there exists an isomorphism of graded algebras
$\eta:\ kG\to \Hom_{A^{\rm e}}(A,A\# kG)$. The composition with the
map $p$ introduced in \leref{3.2a} is an algebra anti-homomorphism
$p\circ \eta:\ kG\to M(A)$. Therefore the map
$$f:\ G\to U(M(A)),~~f(g)=p(\eta(g^{-1}))=: (\rho_{1,g},\rho_{2,g})$$
is a group homomorphism. According to the computations in the proof of
\leref{3.2a}, and using the fact that $\eta(g)$ is homogenous of degree $g$,
we have
$$
f(g)a=\rho_{1,g}(a)=\eta(g^{-1})(g\cdot a)~~{\rm and}~~
af(g^{-1})= \rho_{2,g^{-1}}(a)=\eta(g),$$
and it follows that
$$f(g)af(g^{-1})=(\eta(g)\circ \eta(g^{-1}))(g\cdot a)=g\cdot a,$$
as needed.
\end{proof}

\begin{lemma}\lelabel{3.7}
Let $A$ and $B$ be $G$-module Taylor-Azumaya algebras. Then we have
an isomorphism of graded algebras
\begin{eqnarray*}
&&\hspace*{-2cm}\alpha:\ 
\Hom_{A^{\rm e}}(A,A\# kG)\square \Hom_{B^{\rm e}}(B,B\# kG)\\
&\to & \Hom_{(A\ot B)^{\rm e}}((A\ot B),(A\ot B)\# kG).
\end{eqnarray*}
\end{lemma}

\begin{proof}
The part of degree $g$ of the cotensor product
on the left hand side is generated by $f\# \tilde{f}$, with
$f:\ A\to A\# kG$ and $\tilde{f}:\ B\to B\# kG$ homogeneous of degree
$g$. So we have $f(a)=f_g(a)\# g$, $\tilde{f}(b)=\tilde{f}_g(b)\# g$,
for all $a\in A$ and $b\in B$
(see \leref{3.1a}). Now we define $\alpha(f\ot \tilde{f})$ by
$$\alpha(f\ot \tilde{f})(a\ot b)=(f_g(a)\ot \tilde{f}_g(b))\#g. $$
It is clear that $\alpha(f\ot \tilde{f})$ is a morphism of left $A^{\rm e}$-modules
of degree $g$, and $\alpha$ is a morphism of graded $k$-algebras.
Any morphism of graded $k$-algebras between $kG$-Galois objects is an
isomorphism (see \cite[Prop. 8.1.10]{Ca}), hence $\alpha$ is an isomorphism.
\end{proof}

Collecting the results of all the previous Lemmas, we obtain:

\begin{corollary}\colabel{3.8}
We have a well-defined group homomorphism $\pi:\ \BTM(k,G)\to\Gal(k,G)$, 
given by $\pi[A]=[\Hom_{A^{\rm e}}(A,A\# kG)]$. $\Ker(\pi)$ contains
$\BT(k)$.
\end{corollary}

Let $B$ be a $kG$-Galois object. The inverse map $\beta^{-1}$ of the canonical isomorphism $\beta:~ B\ot B\to B\ot kG$ induces the Miyashita action of $G$ on $B$ (see \cite{Mi}):
\begin{equation}\label{miyashita}
g\rhd b=\sum X^g_i bY^g_i, \ \ \ \sum X^g_i \ot Y^g_i=\beta^{-1}(1\ot g).
\end{equation}
The $G$-graded module $B$ with the Miyashita $G$-action (\ref{miyashita}) is  a Yetter-Drinfeld module, that is, the $G$-action and the $G$-grading satisfy the compatibility relation $g\rhd b_{\sigma}\in B_{g\sigma g^{-1}}$, for every homogeneous element $b_{\sigma}\in B_\sigma$ and $g\in G$. Moreover, $B$ is quantum-commutative in the sense that
\begin{equation}\label{quantumcom}
b_{\sigma} a=(\sigma\rhd a)b_{\sigma}
\end{equation}
for all $a\in B$, $b_{\sigma}\in B_{\sigma}$. 

\begin{lemma}\lelabel{3.8}
The group homomorphism $\pi$ is split surjective.
\end{lemma}

\begin{proof}Let $B$ be a $kG$-Galois object, and consider the multiplier Hopf algebra
$k(G)$ from \exref{1.4}.
We show that $A=B\# k(G)$ is a $G$-Taylor-Azumaya algebra and $\pi(A)\cong B$. Since $A$ 
is clearly flat, $A$ will be equivalent to some finitely generated $G$-Taylor-Azumaya
algebra. By \cite[Thm.4.3]{VDZ}, there exists a strict Morita context
$$\{B\# k(G),k, B, B,[,], (,)\}$$
where $[a,b]=ap_eb=\sum a(p_g\cdot b)\# p_{g^{-1}}$ and $(a,b)=(ab)_e$ for all $a,b\in B$, and $e$ is the identity of $G$. It follows that $A$ is isomorphic to the elementary Taylor-Azumaya algebra $E(\ul{B})$, where $\ul{B}=(B,B,[,])$. So $A$ is a $G$-module Taylor-Azumaya algebra. Define $\phi:\ B\to\pi(A)$ by
$$ \phi(b_{\sigma})=\sum_{g\in G} ((g\rhd b_{\sigma})\# p_g )\# \sigma\in M(A\#kG)$$
for all homogeneous $b_{\sigma}\in B_\sigma$. It is obvious that $\phi$ is graded. 
We show that $\phi$ is a well-defined algebra map.
First we verify that $\phi(B) \subseteq \pi(A)$. 
Let $x=a_{\delta}\# p_h\in A=B\#k(G)$ and $b_{\sigma}\in B_{\sigma}$. 
Identify $x$ with $x\# 1$ in $A\# kG$. Then we have 
\begin{eqnarray*}
&&\hspace*{-1cm}\phi(b_{\sigma})(x\# 1) = \phi(b_{\sigma})[(a_{\delta}\# p_h)\# 1] \\&=&  \sum_{g\in G} [((g\rhd b_{\sigma})\# p_g)\#\sigma ][(a_{\delta}\# p_h) \# 1] = \sum_{g\in G} [((g\rhd b_{\sigma})\#p_g)(a_{\delta}\# p_{h{\sigma^{-1}}})]\# \sigma\\&=& \sum_{g\in G} [(g\rhd b_{\sigma})({p_g}_{(1)}\cdot a_{\delta})\# {p_g}_{(2)}p_{h{\sigma^{-1}}})]\# \sigma=\sum_{g\in G} [(g\rhd b_{\sigma}) a_{\delta})\# p_{\delta^{-1}g}p_{h{\sigma^{-1}}}]\# \sigma \\&=&[(\delta h\sigma^{-1}\rhd b_{\sigma})a_{\delta}\# p_{h{\sigma^{-1}}}]\#\sigma= [a_{\delta}(h\sigma^{-1}\rhd b_{\sigma})\# p_{h{\sigma^{-1}}}]\#\sigma,
\end{eqnarray*}
by (\ref{quantumcom}). On the other hand, we have $g\rhd b_{\sigma}\in B_{g\sigma g^{-1}}$ and
\begin{eqnarray*}
&&\hspace*{-1cm}(x\# 1)\phi(b_{\sigma}) = \sum_{g\in G} [(a_{\delta}\# p_h)\# 1][((g\rhd b_{\sigma})\# p_g )\# \sigma]\\
&=&\sum_{g\in G} [(a_{\delta}\# p_h)((g\rhd b_{\sigma})\# p_g )]\# \sigma
=\sum_{g\in G} [a_{\delta}({p_h}_{(1)}\cdot(g\rhd b_{\sigma}))\# {p_h}_{(2)}p_g )]\# \sigma\\
&=& \sum_{g\in G} [a_{\delta}(g\rhd b_{\sigma})\# p_{g\sigma^{-1}g^{-1}h}p_g )]\# \sigma
= [a_{\delta}(h\sigma^{-1}\rhd b_{\sigma})\# p_{h\sigma^{-1}}]\# \sigma.
\end{eqnarray*}
Thus we have $(x\#1)\phi(b)=\phi(b)(x\# 1)$ for all $x\in A$ and $b\in B$. So $\phi(B)\subseteq \pi(A)$. Next we show that $\phi$ is an algebra map. Let ${a_{\sigma}}\in B_\sigma$ and ${b_{\delta}}\in B_\delta$. Then
\begin{eqnarray*}
&&\hspace*{-1cm}\phi({a_{\sigma}})\phi({b_{\delta}}) = [\sum_{g\in G} ((g\rhd {a_{\sigma}})\# p_g)\#\sigma][\sum_{h\in G} ((h\rhd {b_{\delta}})\# p_h)\# \delta]\\
&=& \sum_{g,h\in G} [(g\rhd {a_{\sigma}})\# p_g][(h\rhd {b_{\delta}})\# p_{h\sigma^{-1}}]\# \sigma\delta\\
&=&\sum_{g,h\in G} [(g\rhd {a_{\sigma}})({p_g}_{(1)}\cdot(h\rhd {b_{\delta}}))\# {p_g}_{(2)}p_{h\sigma^{-1}}]\# \sigma\delta\\
&=& \sum_{g,h\in G} [(g\rhd {a_{\sigma}})(h\rhd{b_{\delta}})\# p_{(h\delta h^{-1})^{-1}g}p_{h\sigma^{-1}}] \# \sigma\delta\\
&=& \sum_{h\in G} [(h\delta\sigma^{-1}\rhd {a_{\sigma}})(h\rhd{b_{\delta}})\#p_{h\sigma^{-1}}]\# \sigma\delta\\
&=& \sum_{t\in G} [(t\delta^{\sigma}\rhd {a_{\sigma}})(t\sigma\rhd{b_{\delta}})\# p_t]\#\sigma\delta=
\sum_{t\in G} [t\rhd((\delta^{\sigma}\rhd {a_{\sigma}})(\sigma\rhd{b_{\delta}}))\# p_t]\#\sigma\delta\\
&=& \sum_{t\in G} [t\rhd((\sigma\rhd{b_{\delta}}){a_{\sigma}})\# p_t] \#\sigma\delta=\sum_{t\in G} [t\rhd({a_{\sigma}}{b_{\delta}})\# p_t]\#\sigma\delta =\phi({a_{\sigma}}{b_{\delta}}),
\end{eqnarray*}
where $\delta^{\sigma}=\sigma\delta\sigma^{-1}$. In the eighth and ninth equality,we used the quantum-commutativity (\ref{quantumcom}) of $B$.Finally, $\phi$ is an isomorphism since $B$ and $\pi(A)$ are $kG$-Galois objects.
\end{proof}

We can now prove the main result of this paper.
\begin{theorem}\thlabel{3.9}
Let $G$ be a group and $k$ a commutative ring. Then we have the following split exact sequence
\begin{equation}\label{main}
1\to \Br'(k)\to \BM'(k,G)\rTo^{\pi} \Gal(k,G)\to 1.
\end{equation}
\end{theorem}
\begin{proof} By \leref{3.8}, it suffices to show that $\Ker(\pi)=\BT(k)$. Let $\iota:~\BT(k)\to \BTM(k,G)$ be the canonical embedding. We define a group morphism 
$\zeta:~\Ker(\pi)\to \BT(k)$ by forgetting the $G$-action on $G$-module 
Taylor-Azumaya algebras. It is easy to see that $\zeta\circ\iota$ is the identity
automorphism of $\BT(k)$. Thus $\zeta$ is a surjective map. We show that $\zeta$ 
is injective as well. Let $A$ be a $G$-module Taylor-Azumaya algebra such that 
the underlying algebra $A$ is elementary, say $A\cong E(\ul{M})$ for some dual pair
of $k$-modules $\ul{M}=(M,M',\mu)$. Since $\pi(A)\cong kG$, it follows from
\leref{3.6} that the $G$-actions on $A$ and on $E(\ul{M})$ are strongly inner. 
By \leref{3.5}, there exists a $G$-module structure on $\ul{M}$ that induces 
the $G$-module structure on $E(\ul{M})$,hence $[E(\ul{M})]=1$.
\end{proof}

If $k$ is a field, $\Br'(k)=\Br(k)$, $\Gal(k,G)=H^2(G,k^*)$, and $\BM'(k,G)=\BM(k,G)$,
and the exact sequence (\ref{main})  is isomorphic to Fr\"ohlich's exact sequence  
\cite[Theorem 4.1]{F}. If
$G$ is finite, then the Picco-Platzeck exact sequence (cf. \cite{PP}) is a `subsequence' of (\ref{main}). 
In this case, $\pi$ restricts to a surjective map on $\BM(k,G)$ and the 
kernel of the restricted map $\pi$ is $\Br(k)$.\\
 
It is not hard to generalize the exact sequence (\ref{main}) to the equivariant Brauer group of an infinite cocommutative Hopf algebra with an integral, which will
yield the infinite version of Beattie's exact sequence \cite{Be}.
The generalization to the equivariant Brauer group of an infinite
quantum group seems far from obvious.

\end{document}